\documentclass[12pt]{article}
\begin{document}

\newpage

\begin{center}
{\large\bf AN ARITHMETIC FUNCTION OF TWO VARIABLES\\}

\large P.A.Gustomesov

\end{center}
\hrule
\begin{abstract}

The arithmetic function of two variables, $Ns \,( a , n )$, is
defined. For all positive integers $n$ and non-negative integers $a$

\begin{displaymath}
Ns \, ( a , n )  = \varphi \, ( n ) \, \frac{ \mu \left(
\frac{ n }{ ( a , n ) } \right) }{ \varphi \left( \frac{ n }{ ( a
, n ) } \right) },
\end{displaymath}
where $\varphi$  is the Euler function, $\mu$  is the
M\"{o}bius function and $(a, n)$ is the greatest common
divisor of integers $a$ and $n$.
Some properties of the function
are given along with the formula that is an analog of the
so-called M\"{o}bius' inversion formula. A heuristic statement is
suggested. The generating function for values of function $Ns \,
( a , n )$  and also a new characteristic property  of prime
numbers are corollaries of the statement.
\end{abstract}
\hrule
\vspace{10mm}

We define the arithmetic function of two variables, $Ns ( a , n )$,
as follows: $n$ stands for all positive integers and $a$ stands for
all non-negative integers,

\begin{displaymath}
Ns ( a , n ) = \varphi ( n ) \frac{ \mu \left( \frac{ n }{ ( a ,
n ) } \right) }{ \varphi \left( \frac{ n }{ ( a , n ) } \right)
},
\end{displaymath}
where $\varphi$ is the Euler function, $\mu$ is the M\"{o}bius
function and $( a , n )$  is the greatest common divisor of
integers $a$ and $n$.
Generalized character of the $Ns ( a , n )$ function is obvious,
since
\begin{displaymath}
Ns ( a , n ) = \varphi ( n ),\,\,\, \mbox{if} \,\,\, ( a , n )
= n ;\,\,\,\,\, Ns ( a , n ) = \mu ( n ) , \,\,\, \mbox{if}\,
\,\, ( a , n ) = 1.
\end{displaymath}
Therefore, one can expect that the function $Ns ( a , n )$
possesses a number of properties, which are both similar to those
of the Euler function and of the M\"{o}bius function as well.

Theorem 1.

\begin{displaymath}
\sum \limits_{d \mid n} Ns ( a , d ) = \left\{ \begin{array}{ll}
n\, , &\,\,\,\,\, \mbox{if $n \mid a$ },\\
0\, , &\,\,\,\,\, \mbox{otherwise},\\
\end{array}
\right.
\end{displaymath}
where the sum is extended to the all divisors $d$ of integer $n$.

Proof. Let the integer $n$ be represented as $ n = n_1 \cdot ( a
, n )$. Then,

\begin{displaymath}
\sum \limits_{d \mid n} Ns ( a , d ) = \sum
\limits_{
d = d_1 \cdot d_2;\,\, d_2 = ( a , d );\,\, d \mid n  } Ns ( a ,
 d_1 \cdot d_2 ) = \sum \limits_{ d_1 \mid n_1;\,\,d_2 \mid ( a , n
 );\,\, d_2 = ( a , d_1 \cdot d_2 ) } \varphi ( d_1 \cdot d_2 ) \frac{ \mu
 ( d_1 ) }{ \varphi ( d_1 ) } =
\end{displaymath}
\begin{displaymath}
 = \sum \limits_{d_1
 \mid n_1} \mu ( d_1 ) \sum \limits_{
 d_2 = d_2'\cdot d_2'';\,\, d_2' = (d_1 , d_2 ) =  ( d_1,\, ( a ,
 n ) );\,\, d_2''  \mid  \frac{ ( a , n ) }{d_2' }
  } \frac{
 \varphi ( d_1 \cdot d_2' \cdot d_2'' ) }{ \varphi ( d_1 ) } = \sum
 \limits_{d_1 \mid n_1} \mu ( d_1 ) \sum \limits_{ d_2'' \mid
 \frac{ ( a , n ) }{ d_2' };\, d_2' = ( d_1, ( a , n ) ) } d_2'
 \varphi ( d_2'' ) = 
\end{displaymath}

\begin{eqnarray}
= \sum \limits_{d_1 \mid n_1 } \mu ( d_1 ) ( a
 , n )
 = \left\{ \begin{array}{ll} ( a , n ) = n
 ,\,\,\, & \mbox{if $n_1 = 1$, i.e.  $n \mid a$};\\\nonumber
 0 ,\,\,\,\ & \mbox{if $n_1 > 1$, i.e. $n \dagger a$}.
 \end{array} \right. &~\\\nonumber
 \end{eqnarray}

Incidently, it follows from the above Proof that
\begin{displaymath}
\sum \limits_{d \mid n} \mid Ns ( a , d ) \mid = ( a , n ) \cdot 2^k\,,
\end{displaymath}
where $k$  is the number of prime divisors of integer $\frac{ n }{ ( a , n )
}$.

The statements similar to the so-called M\"{o}bius' inversion
formulae hold for the function $Ns ( a , n )$,  \cite{1}.

Theorem 2. ( An analog of the first M\"{o}bius' inversion formula).\\
Let $f$ is the arithmetic function and

\begin{displaymath}
g ( n ) = \sum \limits_{d \mid n }  f ( d ) .
\end{displaymath}
Then,
\begin{displaymath}
\sum \limits_{ d \mid ( a , n ) }  f \left( \frac{ n }{ d } \right) \,d =
\sum \limits_{ d \mid n } Ns ( a , d ) \, g \left( \frac{ n }{ d }
\right) .
\end{displaymath}

Proof. We have
\begin{displaymath}
\sum \limits_{d \mid n} Ns ( a , d ) g \left( \frac{ n }{ d } \right) =
\sum \limits_{d \mid n} Ns ( a , d ) \sum \limits_{d' \mid \frac{ n }{ d } } f
(
d' ) = \sum \limits_{d' \mid n} f ( d' ) \sum \limits_{d \mid \frac{ n }{ d' }
}
Ns ( a , d )
\end{displaymath}
and, therefore, according to Theorem 1
\begin{displaymath}
\sum \limits_{ d \mid n }  Ns ( a , d ) \, g \, \left( \frac{ n
}{ d } \right) = \sum \limits_{ \frac{ n }{ d' } \mid a }  f (
d')  \frac{ n }{  d' } = \sum \limits_{ d \mid ( a , n ) } \, f
\, \left( \frac{ n }{ d } \right) \, d.
\end{displaymath}

Corollary 1. ( $ ( a , n ) = 1 $. It is the first M\"{o}bius'
inversion formula).
\begin{displaymath}
f \, ( n ) = \sum \limits_{ d \mid n } \mu ( d ) \, g \,
\left( \frac{ n }{ d } \right).
\end{displaymath}

Corollary 2. ( $( a , n ) = n$ ; $ f ( d ) = d $ ).
\begin{displaymath}
n \,\, d \, ( n ) = \sum \limits_{d \mid n }  \varphi \, ( d ) \,
S \, \left( \frac{  n }{ d } \right),
\end{displaymath}
where $d \, ( n )$  is the number of positive divisors of
integer $n$ and $S$ is the sum of the positive divisors.

Corollary 3. ( $( a , n ) = n\,\,; \,\, f \, (d) = 1$).
\begin{displaymath}
S \, ( n ) = \sum \limits_{ d \mid n } \varphi \, ( d ) \, d \,
\left( \frac{ n }{ d } \right).
\end{displaymath}
And so on.

In general, many expressions containing the M\"{o}bius
function or (and) the Euler function have analogs for the
function $Ns$ as well. For instance, the statement (according to
Theorem 1)

\begin{displaymath}
\sum \limits_{n=1}^{\infty}  \frac{ 1 }{ n^2 }  \cdot \sum
\limits_{m=1}^{ \infty }  \frac{ Ns ( a , m ) }{ m^2 }  = \sum
\limits_{ k \mid a } \frac{ 1 }{ k } ;
\end{displaymath}
or

\begin{displaymath}
 \sum \limits_{m=1}^{ \infty } \frac{ Ns ( a , m ) }{ m^2 }  =
\frac{ 6 }{ \pi^2 } S_{ - 1 } ( a ) ,\,\,\,\,\,\mbox{where}\,\,\,\,\,
S_{ - 1 } ( a ) \equiv \sum \limits_{ d \mid a }  \frac{ 1 }{ d},
\end{displaymath}
is an analog of the statement \cite{1}

\begin{displaymath}
\sum \limits_{n=1}^{\infty}  \frac{ 1 }{ n^2 }\cdot \sum \limits_{
m=1 }^{\infty}  \frac{ \mu (m) }{ m^2 } = \sum \limits_{k=1}^{
\infty }  \frac{ c_k }{ k^2 } = 1,
\end{displaymath}
where $c_k = \sum \limits_{ \ell \mid k } \mu ( \ell )$.

As for summation with respect to $a$, we note the following
property.

Theorem 3. Let $n = p_1^{ \alpha_1 } ... p_k^{ \alpha_k }$
is the canonical expansion. Then,

\begin{displaymath}
\sum \limits_{a=1}^{n}  ( Ns ( a , n ) )^m = ( \varphi ( n ) )^m
\prod \limits_{i=1}^{k}  \left( 1 - \frac{ 1 }{ ( 1 - p_i )^{ m -
1 } } \right),
\end{displaymath}
where $m$  is the non-negative integer.

Proof. It is sufficient to show that

\begin{displaymath}
\sum \limits_{a=1}^{n}  \left( \frac{ \mu \left( \frac{ n }{ ( a
, n ) } \right) }{ \varphi \left( \frac{ n }{ ( a , n ) } \right)
} \right)^m = \left( 1 - \frac{ 1 }{ ( 1 - p )^{ m - 1 } }
\right)  \sum \limits_{ a = 1 }^{n_1} \left( \frac{ \mu
\left( \frac{ n_1 }{ ( a , n_1 ) } \right) }{ \varphi \left(
\frac{ n_1 }{ ( a , n_1 ) } \right) } \right)^m,
\end{displaymath}

where $n = n_1 p^{ \alpha }$; $( n_1, p ) = 1$ ($p$  is the
prime number, $\alpha \geq 1$). We have

\begin{displaymath}
\sum \limits_{a=1}^{n}  \left( \frac{ \mu \left( \frac{ n }{ ( a
, n ) } \right) }{ \varphi \left( \frac{ n }{ ( a , n ) } \right)
} \right)^m =
\sum \limits_{1 \leq a \leq n;\,\,\,(a,p^\alpha) \geq p^{ \alpha
- 1 } }  \left( \frac{ \mu \left( \frac{ n }{ ( a
, n ) } \right) }{ \varphi \left( \frac{ n }{ ( a , n ) } \right)
} \right)^m =
\sum \limits_{a=1}^{n_1 \cdot p }  \left( \frac{ \mu \left( \frac{
n_1 \cdot p }{ ( a
, n_1 \cdot p  ) } \right) }{ \varphi \left( \frac{ n_1 \cdot p  }{ ( a
, n_1 \cdot p  ) } \right)
} \right)^m =
\end{displaymath}
\begin{displaymath}
=\sum \limits_{1 \leq a < n_1\cdot p \, ; \,\, ( a , p ) = 1}  \left(
\frac{ \mu \left( \frac{ n_1 \cdot p }{ ( a
, n_1 \cdot p ) } \right) }{ \varphi \left( \frac{ n_1 \cdot p }{
( a , n_1 \cdot p ) } \right)
} \right)^m
+ \sum \limits_{a=1}^{n_1}  \left( \frac{ \mu \left( \frac{ n_1 }{ ( a
, n_1 ) } \right) }{ \varphi \left( \frac{ n_1 }{ ( a , n_1 ) } \right)
} \right)^m =
( p - 1 )
\sum \limits_{a=1}^{n_1}  \left( \frac{ \mu \left( \frac{
n_1 \cdot p }{ ( a
, n_1  ) } \right) }{ \varphi \left( \frac{ n_1 \cdot p  }{ ( a
, n_1  ) } \right)
} \right)^m +
\sum \limits_{a=1}^{n_1}  \left( \frac{ \mu \left( \frac{ n_1 }{ ( a
, n_1 ) } \right) }{ \varphi \left( \frac{ n_1 }{ ( a , n_1 ) } \right)
} \right)^m =
\end{displaymath}
\begin{displaymath}
=\left( 1 - \frac{ 1 }{ ( 1 - p )^{ m - 1 } } \right)
\sum \limits_{a=1}^{n_1}  \left( \frac{ \mu \left( \frac{ n_1 }{ ( a
, n_1 ) } \right) }{ \varphi \left( \frac{ n_1 }{ ( a , n_1 ) } \right)
} \right)^m .
\end{displaymath}

For negative integers $m$, obviously, the formula

\begin{displaymath}
\sum \limits_{a=1}^{n}  \left( \frac{ \varphi ( n )
}{ \varphi \left( \frac{ n }{ ( a , n ) } \right)
} \right)^m  \left( \mu \left( \frac{ n }{ ( a , n ) } \right)
\right)^{ \mid m \mid } = \left( \varphi ( n ) \right)^m \, \prod
\limits_{i=1}^{k} \left( 1 - \frac{ 1 }{ ( 1 - p_i )^{ m - 1 } }
\right)
\end{displaymath}
is valid.

It follows from Theorem 3 that

\begin{displaymath}
\sum \limits_{a=1}^{n} Ns ( a , n ) = \left\{ \begin{array}{ll}
1\,\,\, &, \,\,\,\,\, \mbox{if}\,\,\, n = 1, \\
0\,\,\, &, \,\,\,\,\, \mbox{if}\,\,\, n > 1;
\end{array}
\right.
\end{displaymath}
and also
\begin{displaymath}
\sum \limits_{a=1}^{n}  \mid Ns ( a , n ) \mid = \varphi ( n )
\, 2^k,
\end{displaymath}
where $k$  is the number of prime divisors of the integer $n$.

The following property of the function $Ns (a , n )$ is of
particular interest. It is given here in the form of hypothesis
because the author has no a completed proof at his disposal.

Theorem 4. Let

\begin{displaymath}
\prod \limits_{i=1}^{n-1}  ( 1 - q^i ) = \sum \limits_{k=0}^{
\frac{ n ( n - 1 ) }{2} }  p_{ n - 1 } ( k ) \, q^k.
\end{displaymath}
(Coefficients $p_{ n - 1 } ( k )$ are of a specified sense in the
theory of partitions  \cite{2}). Then,\\
for $0 \leq a \leq ( n - 1 )$

\begin{displaymath}
Ns ( a , n ) = \sum \limits_{ k \geq 0;\,\,\, a + n\cdot k \leq
\frac{ n ( n - 1 ) }{ 2 } }  p_{ n - 1 } ( a + n \cdot k ) = \sum
\limits_{k=0}^{ \left[ \frac{ n - 1 }{ 2 } - \frac{ a }{ n }
\right] }  p_{ n - 1 } ( a + n \cdot k ),
\end{displaymath}
where $\left[ \frac{ n - 1 }{ 2 } - \frac{ a }{ n } \right]$
is the integral part of the number.

Corollary 1.

\begin{displaymath}
\varphi ( n ) = \sum \limits_{k = 0}^{ \left[ \frac{ n - 1 }{ 2 }
\right] }  p_{ n - 1 } ( n \cdot k );
\end{displaymath}
\begin{displaymath}
\mu ( n ) = \sum \limits_{k = 0}^{ \left[ \frac{ n - 2 }{ 2 }
\right] }  p_{ n - 1 } ( 1 + n \cdot k ).
\end{displaymath}

Corollary 2. (It expresses a characteristic property of prime
numbers.) If $ p $  is a prime odd number, then

\begin{displaymath}
1 + \sum \limits_{k = 0}^{ \left( \frac{ p - 1 }{ 2 }
\right) }  p_{ n - 1 } ( p \cdot k ) = p,
\end{displaymath}
and also
\begin{displaymath}
1 + \sum \limits_{k = 0}^{ \left( \frac{ p - 1 }{ 2 } - 1 \right) }
p_{ n - 1 } ( a + p \cdot k ) = 0,
\end{displaymath}
where $1 \leq a \leq ( p - 1 )$.

Corollary 3. (It is a generating function for values of the
function  $Ns ( a ,n )$.) Let

\begin{displaymath}
\frac{ ( - 1 )^{ n - 1 }}{ 1 - q^n } \prod \limits_{i=1}^{n-1}  (1
- q^i ) = \sum \limits_{k=0}^{ \infty }  N_n ( k ) \, q^k.
\end{displaymath}
Then, for $a \geq 1$

\begin{displaymath}
Ns ( a , n ) = N_n \left( \frac{ ( n - 1 ) \, ( n - 2 ) }{ 2 } +
a - 1 \right) .
\end{displaymath}

\end{document}